\documentclass[8 pt]{amsproc}
\usepackage{amsmath}
\usepackage{graphicx}
\usepackage{amssymb}
\usepackage[active]{srcltx}
\oddsidemargin=-1cm\evensidemargin=-1cm
\begin{document}
\numberwithin{equation}{section}

\def\1#1{\overline{#1}}
\def\2#1{\widetilde{#1}}
\def\3#1{\widehat{#1}}
\def\4#1{\mathbb{#1}}
\def\5#1{\frak{#1}}
\def\6#1{{\mathcal{#1}}}

\def\C{{\4C}}
\def\R{{\4R}}
\def\n{{\4n}}
\def\Z{{\4Z}}

\title{Formal   Equivalences in $\mathbb{C}^{4}$}
\author{Valentin Burcea}
\begin{abstract} There are  considered   formal constructions of normal form type using Formal (Holomorphic Segre-Preserving) Equivalences, for complexifications of Real-Formal Hypersurfaces in $\mathbb{C}^{2}$. Such formal constructions   are convergent   if the source manifolds are Analytic.   
\end{abstract}
\address{V. Burcea: INDENDENT  }
\email{vdburcea@gmail.com}

\thanks{\emph{Keywords:} Finite Jet Determination,  Normal Form, Cauchy-Riemann Geometry, CR Equivalence, Equivalence Problem}
\thanks{Special Thanks (in regard to  paper)  to  Science Foundation Ireland grant 10/RFP/MT H2878}
\thanks{I make clear that the reference \cite{V1} was fully supported by  Science Foundation Ireland Grant 06/RFP/MAT 018}

\maketitle
 
\def\Label#1{\label{#1}{\bf (#1)}~}


\def\cn{{\C^n}}
\def\cnn{{\C^{n'}}}
\def\ocn{\2{\C^n}}
\def\ocnn{\2{\C^{n'}}}


\def\dist{{\rm dist}}
\def\const{{\rm const}}
\def\rk{{\rm rank\,}}
\def\id{{\sf id}}
\def\tr{{\bf tr\,}}
\def\aut{{\sf aut}}
\def\Aut{{\sf Aut}}
\def\CR{{\rm CR}}
\def\GL{{\sf GL}}
\def\Re{{\sf Re}\,}
\def\Im{{\sf Im}\,}
\def\span{\text{\rm span}}
\def\Diff{{\sf Diff}}

\def\codim{{\rm codim}}
\def\crd{\dim_{{\rm CR}}}
\def\crc{{\rm codim_{CR}}}

\def\phi{\varphi}
\def\eps{\varepsilon}
\def\d{\partial}
\def\a{\alpha}
\def\b{\beta}
\def\g{\gamma}
\def\G{\Gamma}
\def\D{\Delta}
\def\Om{\Omega}
\def\k{\kappa}
\def\l{\lambda}
\def\L{\Lambda}
\def\z{{\bar z}}
\def\w{{\bar w}}
\def\Z{{\1Z}}
\def\t{\tau}
\def\th{\theta}

\emergencystretch15pt \frenchspacing

\newtheorem{Thm}{Theorem}[section]
\newtheorem{Cor}[Thm]{Corollary}
\newtheorem{Pro}[Thm]{Proposition}
\newtheorem{Lem}[Thm]{Lemma}

\theoremstyle{definition}\newtheorem{Def}[Thm]{Definition}

\theoremstyle{remark}
\newtheorem{Rem}[Thm]{Remark}
\newtheorem{Exa}[Thm]{Example}
\newtheorem{Exs}[Thm]{Examples}

\def\bl{\begin{Lem}}
\def\el{\end{Lem}}
\def\bp{\begin{Pro}}
\def\ep{\end{Pro}}
\def\bt{\begin{Thm}}
\def\et{\end{Thm}}
\def\bc{\begin{Cor}}
\def\ec{\end{Cor}}
\def\bd{\begin{Def}}
\def\ed{\end{Def}}
\def\br{\begin{Rem}}
\def\er{\end{Rem}}
\def\be{\begin{Exa}}
\def\ee{\end{Exa}}
\def\bpf{\begin{proof}}
\def\epf{\end{proof}}
\def\ben{\begin{enumerate}}
\def\een{\end{enumerate}}
\def\beq{\begin{equation}}
\def\eeq{\end{equation}}

\section{Introduction and Main Results}
 
Let $M,\hspace{0.1 cm}N\subset \mathbb{C}^{2}$ two Real-Formal Hypersurfaces   defined near $0\in\mathbb{C}^{2}$ by
\begin{equation}M:\rho_{1}\left(Z,\overline{Z}\right)=0,\quad N:\rho_{2}\left(Z,\overline{Z}\right)=0,\label{V11}
\end{equation}
where $Z=\left(w,z\right)$ are coordinates in $\mathbb{C}^{2}$.

We recall from Angle\cite{An} that any (Formal) Holomorphic Segre Preserving Mapping, between $M$ and $N$, is defined like
\begin{equation}\mathcal{H}:\mathbb{C}^{2}\rightarrow \mathbb{C}^{2}\quad\mbox{such that $\mathcal{H}\left(Z,\zeta\right)=\left(H\left(Z\right),\tilde{H}\left(\zeta\right)\right)$,} \label{cete}
\end{equation} 
 where $H,\tilde{H}:\mathbb{C}^{2}\rightarrow \mathbb{C}^{2}$ are formal holomorphic mappings such that  
\begin{equation}\rho_{2}\left(H\left(Z\right),\tilde{H}\left(\zeta\right) \right)=0,\quad\mbox{ for all $\left(Z,\zeta\right)\in \mathbb{C}^{4}$ such that $\rho_{1}\left(Z,\zeta\right)=0$,}\label{cete11}
\end{equation}
where $\zeta=\left(\nu,\xi\right)$ replaces $Z=\left(w,z\right)$.

Two   (Formal) Holomorphic Segre Preserving Mappings $\mathcal{H}_{1},\hspace{0.1 cm}\mathcal{H}_{2}$     are called to be determined by  $k$-jets if   
$ J^{k}\left(\mathcal{H}_{1}\right)=J^{k}\left(\mathcal{H}_{2}\right)\rightarrow \mathcal{H}_{1}=\mathcal{H}_{2}$, 
where $k\in\mathbb{N}^{\star}$ and $J^{k}$ defines the $k$-jet in the corresponding formal expansion.

This paper studies standard questions from \cite{Lam}  in order to understand  the Formal Holomorphic Segre Preserving Mappings\cite{An},\cite{Zha} of non-trivial Real-Formal Hypersurfaces in   $\mathbb{C}^{2}$. In particular, we use Formal  (Holomorphic)  Segre Equivalences in order to consider (formally) constructions of normal form type recalling the  methods learned by the author\cite{V1},\cite{V2}   from Zaitsev\cite{D1},\cite{D3}. In particular, we obtain:
\bt\label{t} Let $M\subset\mathbb{C}^{2}$ be a Real-Formal Hypersurface defined near $p=0$ by
\begin{equation}
  \Im w=\left(\Re w\right)^{s}\mathcal{L}\left(z,\overline{z}\right)+\displaystyle\sum_{m+n+p\geq k_{0}+1 }\varphi_{mnp}{z}^{m}{\overline{z}}^{n}\left(\Re w \right)^{p},\label{vvv}
\end{equation}
where $\mathcal{L}\left(z,\overline{z}\right)$ is a   homogeneous polynomial of degree $k_{0}-s\geq 3$, for given $k_{0}, s\in\mathbb{N}$  such that 
 $\varphi_{mn0}=0$ if $s\neq 0$, for all $m,n\in\mathbb{N}$ with   $m+n\geq k_{0}+1$. 
Then,   the Formal Holomorphic Segre Equivalences of $M$ are determined by  their $1$-jets.  
\et 

This case  (\ref{vvv}) defines an  infinite type hypersurface\cite{Mey},\cite{LJu} for $s\neq 0$, and respectively  a finite type hypersurface\cite{K1},\cite{K2} for $s=0$. Normalizations, derived from Weighted and Pseudo-Weighted Versions of the (Generalized) Fischer Decomposition\cite{Sh}, are   formally imposed in  complexified (formal)  local defining equations. Such  Decompositions are applied with respect to a natural System of Weights when $s=0$. Otherwise, we introduce by \cite{V2} a System of Pseudo-Weights   defined   on individual classes of terms.  It verifies   that the (pseudo-)weight of the sum of two terms is actually the sum of their (pseudo-)weights, but the (pseudo-)weight of a product of two terms may not the sum of their (pseudo-)weights. Normalizations  on sums of weighted and (pseudo-)weighted terms are applied using the strategy from \cite{V2},   after a preliminary normalization   of the  $1$-jets of the Formal  Segre Preserving Change of Coordinates. In particular, we obtain:
\bc\label{c1} Any Formal Segre-Preserving Mapping, between  two non-flat Real-Formal Hypersurfaces in $\mathbb{C}^{2}$,  posses the finite jet determination property. 
\ec

The finite jet determination follows as in the standard cases from Ebenfelt-Lamel-Zaitsev\cite{ELZ2},\cite{ELZ3}, regardless of such non-trivialities existent in both cases $s=0$ and $s\neq 0$. It defines a specific  property occuring  in the   theory of    Real Submanifolds in Complex Spaces\cite{BERbook}, which       may not generally  hold. The example of the  group of biholomorphisms  of   Levi-flat hypersurfaces  shows the nontriviality of the finite jet determination problem\cite{ELZ1} in $\mathbb{C}^{2}$ (see also \cite{KLS}).   It is indicated  \cite{BERbook} for an extended introduction to this topic\cite{ELZ2},\cite{ELZ3},\cite{Lam}, and Angle\cite{An},\cite{An1} and Zhang\cite{Zha} for  progresses concerning the Segre Holomorphic Preserving Mappings\cite{An} and related to standard problems in CR Geometry\cite{BERbook}.  

Going forward, it is studied the problem of convergence\cite{BEL3}  of a Formal Equivalence\cite{D2}. In particular, we obtain: 
\bt\label{t2} Any Formal Segre-Equivalence, of two non-flat  Real-Analytic Hypersurfaces in $\mathbb{C}^{2}$, is Convergent. 
\et

It  may be an evidence of the inflexibility of such Formal Equivalences, regardless of the geometrical context.  Such conclusion contrasts when  $s\neq 0$ in   the standard case in $\mathbb{C}^{2}$, where     Kossovskiy-Shafikov\cite{RI} showed  that it may not  exist   Holomorphic Equivalences between two non-minimal Real-Analytic  Formally  Equivalent  Submanifolds in   Complex Spaces. In particular, we obtain
\bc\label{c2} Any Formal Segre-Preserving Mapping, of two non-flat Real-Analytic Hypersurfaces in $\mathbb{C}^{2}$, is Convergent. 
\ec

The proof of Theorem \ref{t2} is simple. It is obtained from an analytic equation derived from the local defining equation.   Moreover, the formal  constructions of normal form type are   convergent   when the source manifolds are analytic (see also \cite{CM},\cite{KZ}). It  is surprising, because Kol\`{a}\v{r}\cite{K2} constructed an example of divergent normal form  for a Real-Analytic Hypersurface\cite{K1} of finite type in $\mathbb{C}^{2}$.  

\subsection{Acknowledgements} It is my own work effectuated independently and derived from my doctoral efforts.  I owe  to Science Foundation of Ireland, because I was using Irish Funding during my doctoral studies in Trinity College Dublin. Special Thanks  to  my supervisor (Prof. Dmitri Zaitsev) for many  conversations regarding  the main part\cite{V1} (fully supported by Science Foundation Ireland Grant 06/RFP/MAT 018) of my doctoral thesis. Special Thanks  also to   Jasmin Raissy, George Ioni\c{t}a, Diogo Bessam and Florian Bertrand.   

\section{Preliminaries}
\subsection{Local Coordinates} It is used the standard  procedure\cite{V1},\cite{V2}  implemented for constructing (formal)  normal forms\cite{K1},\cite{D1},\cite{D3}  in Complex Analysis. In particular, we  work by (\ref{vvv}) in the following local coordinates:

Let $M \subset\mathbb{C}^{2}$ be the  Real-Formal Hypersurface\cite{BERbook}, defined  by
\begin{equation} M:\quad \Im w=\left(\Re w\right)^{s}\mathcal{L}\left(z,\overline{z}\right)+\displaystyle\sum_{m+n+p\geq k_{0}+1}\varphi_{mnp}{z}^{m}{\overline{z}}^{n}\left(\Re w \right)^{p},  \label{ec}
\end{equation} 
where  $\left(w,z\right)$ are coordinates in $\mathbb{C}^{2}$.

Let $M'\subset\mathbb{C}^{2}$ be another  Real-Formal Hypersurface\cite{BERbook}, defined  by
\begin{equation}\quad\quad M':\quad   \Im w'=\left(\Re w' \right)^{s}\mathcal{L}\left(z',\overline{z'}\right)+\displaystyle\sum_{m+n+p\geq k_{0}+1 }{\varphi'}_{mnp}{z'}^{m}{\overline{z'}}^{n}\left(\Re w'\right)^{p},  \label{ec1}
\end{equation} 
where   $\left(w',z'\right)$ are coordinates in $\mathbb{C}^{2}$.

 In order to study   formal (holomorphic) equivalences like (\ref{cete}), we complexify (\ref{ec}) and (\ref{ec1}).  We obtain
\begin{equation}  \mathcal{M}: \quad \frac{w-\nu}{2\sqrt{-1}}=\left(\frac{w+\nu}{2}\right)^{s}\mathcal{L}\left(z,\xi\right)+\displaystyle\sum_{m+n+p\geq k_{0}+1 }\varphi_{mnp}{z}^{m}{\xi}^{n}\left(\frac{w+\nu}{2} \right)^{p},  \label{gigi}
\end{equation}
and respectively, we obtain
\begin{equation}  \quad\quad\hspace{0.2 cm} \mathcal{M}':\hspace{0.1 cm} \frac{ w'-\nu'}{2\sqrt{-1}}=\left(\frac{w'+\nu'}{2}\right)^{s}\mathcal{L}\left(z',\xi'\right)+\displaystyle\sum_{m+n+p\geq k_{0}+1 }\varphi'_{mnp}{z'}^{m}{\xi'}^{n}\left(\frac{w'+\nu'}{2} \right)^{p},  \label{gigi1}
\end{equation}
where   $\overline{z}'$ has been replaced with $\xi'$, and   $\overline{w}'$ has been replaced with $\nu'$. 

The Theorem of Implicit Functions applied in (\ref{gigi}) gives  $w=Q\left(z,\xi,\nu\right)$, 
 where $Q$ is (formal) holomorphic in $\left(z,\xi,\nu\right)$   written   like
\begin{equation}Q\left(z,\xi,\nu\right)=\nu+Q_{k_{0}}\left(z,\xi,\nu\right)+Q_{k_{0}+1}\left(z,\xi,\nu\right),\label{Qu}
\end{equation}
where we have used a  homogeneous polynomial  of degree $k_{0}$ in $\left(z,\xi,\nu\right)$ and     a (possibly infinite) sum of  homogeneous polynomials  of degree at least $k_{0}+1$ in $\left(z,\xi,\nu\right)$,  denoted by  $Q_{k_{0}}$ and  $Q_{k_{0}+1}\left(z,\xi,\nu\right)$, in order to  replace (\ref{Qu}) in (\ref{gigi}). We obtain
\begin{equation}\begin{split}&   \frac{ Q_{k_{0}}\left(z,\xi,\nu\right)+Q_{k_{0}+1}\left(z,\xi,\nu\right) }{2\sqrt{-1}}=  \left(\frac{2\nu+Q_{k_{0}}\left(z,\xi,\nu\right)+Q_{k_{0}+1}\left(z,\xi,\nu\right)}{2}\right)^{s}\mathcal{L}\left(z,\xi\right)\\& \quad\quad\quad\quad\quad\quad\quad\quad\quad\quad\quad\quad\quad\quad\quad\quad\quad\quad\quad\quad\quad \hspace{0.15 cm}   \begin{tabular}{l} \rotatebox[origin=c]{270}{$+$}\end{tabular} \\&\quad\quad\quad\quad\quad\quad\quad\quad\hspace{0.2 cm} 
 \displaystyle\sum_{m+n+p\geq k_{0}+1 }\varphi_{mnp}{z}^{m}{\xi}^{n}\left(\frac{2\nu+Q_{k_{0}}\left(z,\xi,\nu\right)+Q_{k_{0}+1}\left(z,\xi,\nu\right)}{2}\right)^{p}.\end{split} \label{gigi2}
\end{equation}

Then, we consider formal expansions   in (\ref{gigi2}) according to following evaluations
\begin{equation}\begin{split}& \displaystyle\sum_{m+n+p\geq k_{0}+1 }\varphi_{mnp}{z}^{m}{\xi}^{n}\left(\frac{2\nu+Q_{k_{0}}\left(z,\xi,\nu\right)+Q_{k_{0}+1}\left(z,\xi,\nu\right)}{2} \right)^{p}=\mbox{O}\left(k_{0}+1\right),\\& \hspace{0.27 cm} \quad\quad  \quad \quad \quad \quad \quad \left(\frac{2\nu+Q_{k_{0}}\left(z,\xi,\nu\right)+Q_{k_{0}+1}\left(z,\xi,\nu\right)}{2}\right)^{s}\mathcal{L}\left(z,\xi\right)=\nu^{s}\mathcal{L}\left(z,\xi\right)+\mbox{O}\left(k_{0}+1\right).\end{split}\label{gigi21}
\end{equation}
 
After a study of terms in the both sides in (\ref{gigi2}) using (\ref{gigi21}), we obtain
\begin{equation}Q_{k_{0}}\left(z,\xi,\nu\right)=2\sqrt{-1} \nu^{s} \mathcal{L}\left(z,\xi\right).\label{lolo1}
\end{equation}

Next, in order to assume that  the linear part,  of the formal (holomorphic) equivalence $\mathcal{H}$, is standard, or equivalently by (\ref{cete}), that the formal power series $H$ and $\tilde{H}$ have standard linear parts, we  use:  

\subsection{Linear (holomorphic Segre preserving) Changes of Coordinates} 
We write the formal expansions of $H$ and $\tilde{H}$:
\begin{equation}  H(w,z)=\left(\displaystyle\sum_{ k+l\geq 1}f_{ kl}z^{k} w^{l},\displaystyle\sum_{ k+l\geq 1}g_{kl}z^{k} w^{l}\right),\quad\tilde{H}\left(\nu,\xi\right)=\left( \displaystyle\sum_{ \tilde{k}+\tilde{l}\geq 1}\tilde{f}_{\tilde{k}\tilde{l}} \xi^{\tilde{k}} \nu^{\tilde{l}},\displaystyle\sum_{ \tilde{k}+\tilde{l}\geq 1}\tilde{g}_{\tilde{k}\tilde{l}} \xi^{\tilde{k}} \nu^{\tilde{l}}\right), \label{cete1}
\end{equation}
using the   hypothesis that $\mathcal{H}$ is an equivalence, which implies
$f_{10}\neq 0$,   $g_{01}\neq 0$ and $\tilde{f}_{10}\neq 0$,  $\tilde{g}_{01}\neq 0$. 

Next, we replace (\ref{Qu}) in (\ref{cete1}), and then we replace (\ref{cete1}) in (\ref{gigi1}).   Because there do  not exist   multiplications of $\nu$, $z$ and $\xi$ with complex coefficients in right-hand side from  (\ref{gigi1}),  we obtain $g_{10}=\tilde{g}_{10}$ and $g_{01}=\tilde{g}_{01}=0$. 
Then, in order to normalize $\mathcal{H}$, it has sense to consider the following linear change of coordinates
\begin{equation}\left(w',z';\nu',\xi'\right)=\left(\frac{w}{g_{10}}, \frac{z}{\left(g_{10}\right)^{\frac{s}{k_{0}-s}}};\frac{\nu}{\tilde{g}_{10}},\frac{\xi}{\left(\tilde{g}_{10}\right)^{\frac{s}{k_{0}-s}}}\right),\label{18}
\end{equation}  
which is well-defined and preserves   the Model
\begin{equation}
w=\nu+2\sqrt{-1} \nu^{s} \mathcal{L}\left(z,\xi\right).\label{mood}
\end{equation}
 
 Changing  the coordinates using (\ref{18}), we can assume   $g_{10}=\tilde{g}_{10}=1$. Then,
  we identify the terms of degree $2$ in $\left(z,\xi\right)$  in (\ref{gigi1}).  Then 
$$\mathcal{L}\left(z,\xi\right)=\mathcal{L}\left(f_{10}z,\tilde{f}_{10}\xi\right),\quad\mbox{and then}\hspace{0.1 cm} \mathcal{L}\left(\frac{z}{f_{10}},\frac{\xi}{\tilde{f}_{10}}\right)=\mathcal{L}\left(z,\xi\right).$$

Model (\ref{mood}) is preserved by the following  linear (holomorphic Segre preserving) change of coordinates
\begin{equation}\left(w',z';\nu',\xi'\right)=\left(w,\frac{z}{f_{10}} ; \nu , \frac{\xi}{\tilde{f}_{10}} \right).\label{181}
\end{equation} 
 
Changing the coordinates using (\ref{181}), we  can  assume $f_{10}=\tilde{f}_{10}=1$. Then, 
 replacing (\ref{lolo1}) and (\ref{cete1}) in (\ref{gigi1}), we move forward in order to derive  formal constructions of normal form type. In particular, we  study the following: 
\begin{equation}\begin{split}& \quad\quad \frac{\nu+2\sqrt{-1} \nu^{s} \mathcal{L}\left(z,\xi\right)+Q_{k_{0}+1}\left(z,\xi,\nu\right)-\nu+ \displaystyle\sum_{ k+l\geq 2}\left(g_{kl}z^{k}\left(\nu+2\sqrt{-1} \nu^{s} \mathcal{L}\left(z,\xi\right)+Q_{k_{0}+1}\left(z,\xi,\nu\right)\right)^{l}-\tilde{g}_{kl}\xi^{k}\nu^{l}\right)}{2\sqrt{-1}} \\&  \quad\quad\quad\quad\quad\quad\quad\quad\quad\quad\quad\quad\quad\quad\quad\quad\quad
\quad\quad\quad\quad\quad\quad\quad\quad \begin{tabular}{l} \rotatebox[origin=c]{270}{$=$}\end{tabular} \\& \quad\quad\quad\quad\left(\frac{1}{2}\left( \nu^{s} \mathcal{L}\left(z,\xi\right)+  Q_{k_{0}+1}\left(z,\xi,\nu\right)+2\nu+\displaystyle\sum_{k+l\geq 2}\left(g_{kl}z^{k}\left(\nu+Q_{k_{0}}\left(z,\xi,\nu\right)+Q_{k_{0}+1}\left(z,\xi,\nu\right)\right)^{l}+\tilde{g}_{kl}\xi^{k}\nu^{l}\right)\right)\right)^{s}\cdot \\&\quad\quad\quad\quad\quad\quad\quad\quad\quad\quad\quad\quad \quad\quad\quad \quad\quad\quad   \quad\quad\quad\quad\quad\quad\quad\quad\hspace{0.1 cm}\mathcal{L}\left( \displaystyle\sum_{k+l\geq 1}f_{kl}z^{k}\left(\nu+2\sqrt{-1} \nu^{s} \mathcal{L}\left(z,\xi\right)+Q_{k_{0}+1}\left(z,\xi,\nu\right)\right)^{l}, \right.\\& \left. \displaystyle\sum_{k+l\geq 1}\tilde{f}_{kl}\xi^{k}\nu^{l}\right)+\displaystyle\sum_{m+n+p\geq k_{0}+1}\varphi'_{mnp}\left(   \displaystyle\sum_{k+l\geq 1}  f_{kl}z^{k}\left(\nu+2\sqrt{-1} \nu^{s} \mathcal{L}\left(z,\xi\right)+Q_{k_{0}+1}\left(z,\xi,\nu\right)\right)^{l}\right)^{m}\cdot\left( \displaystyle\sum_{k+l\geq 1}\tilde{f}_{kl}\xi^{k}\nu^{l}  \right)^{n} \cdot \\& \left(\frac{1}{2}\left(2\sqrt{-1} \nu^{s} \mathcal{L}\left(z,\xi\right)+  Q_{k_{0}+1}\left(z,\xi,\nu\right)+2\nu +\displaystyle\sum_{k+l\geq 2}\left(g_{kl}z^{k}\left(\nu+2\sqrt{-1} \nu^{s} \mathcal{L}\left(z,\xi\right)+Q_{k_{0}+1}\left(z,\xi,\nu\right)\right)^{l}+\tilde{g}_{kl}\xi^{k}\nu^{l}\right) \right)\right)^{p}.\end{split}\label{e3FFF}
\end{equation}

In order to determine uniquely (\ref{cete1}) from (\ref{e3FFF}), we use:
\subsection{Fischer Decompositions\cite{Sh}} Any formal (holomorphic) power series, denoted as $F\left(z,\xi,\nu\right)$, may be uniquely written   as 
\begin{equation}F\left(z,\xi,\nu\right)=G\left(z,\xi,\nu\right)P\left(z,\xi,\nu\right)+R\left(z,\xi,\nu\right),\quad\mbox{where $P^{\star}\left(R\left(z,\xi,\nu\right)\right)=0$,}\label{de}
\end{equation}
where $G\left(z,\xi,\nu\right)$ and $R\left(z,\xi,\nu\right)$ are  formal  power series in $\left(z,\xi,\nu\right)$, and 
\begin{equation}P^{\star}:=\displaystyle\sum_{m+n+p=r}\overline{p_{mnp}}\frac{\partial^{r}}{\partial z^{m}\partial \xi^{n}\partial\nu^{p}},\quad\mbox{if}\hspace{0.1 cm}P\left(z,\xi,\nu\right)=\displaystyle\sum_{m+n+p=r }p_{mnp}z^{m}\xi^{n}\nu^{p},\quad\mbox{where $r\in\mathbb{N}^{\star}$}.\label{lele}
\end{equation}

Then,  (\ref{de}) can be extended,  for any homogeneous polynomials $\left(P_{1},P_{2},P_{3},P_{4}\right)\left(z,\xi,\nu\right)$  of   degree $r$, in order to write $F\left(z,\xi,\nu\right)$ as:
\begin{equation}F\left(z,\xi,\nu\right)=\displaystyle\sum_{i=1}^{4}G_{i}\left(z,\xi,\nu\right)P_{i}\left(z,\xi,\nu\right)+R\left(z,\xi,\nu\right),\quad\mbox{where $R\left(z,\xi,\nu\right) \in \displaystyle\bigcap_{i=1}^{4}\ker P^{\star}_{i} $,}\label{degen}
\end{equation}
where $\left(G_{1},G_{2},G_{3},G_{4}\right)\left(z,\xi,\nu\right)$ are  formal  power series in $\left(z,\xi,\nu\right)$. 

Also,  (\ref{de}) can be extended,  for any homogeneous polynomials $\left(P_{1},P_{2},P_{3},P_{4},P_{5},P_{6}\right)\left(z,\xi,\nu\right)$  of   degree $r$, in order to write $F\left(z,\xi,\nu\right)$ as:
\begin{equation}F\left(z,\xi,\nu\right)=\displaystyle\sum_{i=1}^{6}G_{i}\left(z,\xi,\nu\right)P_{i}\left(z,\xi,\nu\right)+R\left(z,\xi,\nu\right),\quad\mbox{where $R\left(z,\xi,\nu\right) \in \displaystyle\bigcap_{i=1}^{6}\ker P^{\star}_{i} $,}\label{degenn}
\end{equation}
where $\left(G_{1},G_{2},G_{3},G_{4},G_{5},G_{6}\right)\left(z,\xi,\nu\right)$ are  formal  power series in $\left(z,\xi,\nu\right)$.

The uniqueness of  $R\left(z,\xi,\nu\right)$ is known in (\ref{degenn}) like in (\ref{de}), but it is not generally clear the uniqueness of  $\left(G_{1},G_{2},G_{3},G_{4}\right)\left(z,\xi,\nu\right)$ and 
 $\left(G_{1},G_{2},G_{3},G_{4},G_{5},G_{6}\right)\left(z,\xi,\nu\right)$   in (\ref{degen}) and (\ref{degenn}). In particular, we study:  

\section{Interactions of terms}

\subsection{\mbox{The case $s=0$}} We study  terms in (\ref{e3FFF}) using   Model (\ref{mood}), when $s=0$. In particular, we use the  Model
 \begin{equation}w=\nu+2\sqrt{-1}    \mathcal{L}\left(z,\xi\right),\label{ref1}
\end{equation}

Procedures from Kol\`{a}\v{r}\cite{K1}, Zaitsev\cite{D1},\cite{D3} and Huang-Yin\cite{HY}
are implemented in order  to consider  normalizations   in  (\ref{e3FFF}), when $s=0$. More precisely, it is required to identify   interactions containing the undetermined terms of  (\ref{cete1}) by  linearising in (\ref{e3FFF}). On the other hand, it is not clear which   normalizations are the most suitable, because there exist many  terms available in order to consider them like the vanishing in (\ref{e3FFF}) of the coefficient  of $\nu^{l}$, which   uniquely determines
\begin{equation}g_{kl}\hspace{0.1 cm}\mbox{and}\hspace{0.1 cm} g_{N,0}-\tilde{g}_{N,0}, \quad \mbox{ for all $l\in\mathbb{N}^{\star}$ and $k\in\mathbb{N}$ with $N=k+l\geq 2$.}\label{lolo11}
\end{equation}

Such normalizations are not  satisfactory, because (\ref{cete})  can not be computed entirely, since  there are left undetermined an infinite number of parameters. Moreover,   there exist better  options  which may be  considered   in order to consider   normalizations in (\ref{e3FFF}). In particular,    the  Model (\ref{ref1}) provides new homogeneous terms, but the  Model (\ref{ref1}) is not homogeneous. In particular,  the  Model (\ref{ref1}) can not be  considered in order   to use   classical Fischer Decompositions\cite{Sh}, but the  Model (\ref{ref1}) becomes homogeneous   using the following  System of Weights. 

We define
$\mbox{wt}\left\{\nu\right\}=k_{0}$ and $\mbox{wt}\left\{z\right\}=\mbox{wt}\left\{\xi\right\}=1$ in order    evaluate the weights of the  interactions of terms from (\ref{e3FFF}). In particular, we focus on the following terms
\begin{equation} \begin{split}&\quad
g_{kl}z^{k}\left(\nu+2\sqrt{-1}   \mathcal{L}\left(z,\xi\right)+Q_{k_{0}+1}\left(z,\xi,\nu\right)\right)^{l},\quad\hspace{0.05 cm}\mbox{for $k,l\in\mathbb{N}$ with $N=k+k_{0}l$,}\\& g_{k'l'}z^{k'}\left(\nu+2\sqrt{-1}  \mathcal{L}\left(z,\xi\right)+Q_{k_{0}+1}\left(z,\xi,\nu\right)\right)^{l'},\quad\mbox{for $k',l'\in\mathbb{N}$ with  $N'=k' + k_{0}l'$ and  $2 \leq  N < N' $.} \end{split} \label{lolo1q}
\end{equation}

These terms  (\ref{lolo1q})  do not overlap because of their different pseudo-weights, in the light of the following  weighted-evaluations
\begin{equation*}\begin{split}&\hspace{0.1 cm}   \left(\nu+2\sqrt{-1}   \mathcal{L}\left(z,\xi\right)+Q_{k_{0}+1}\left(z,\xi,\nu\right)\right)^{l}=\left(\nu+2\sqrt{-1}  \mathcal{L}\left(z,\xi\right)\right)^{l}+\mbox{O}_{\mbox{wt} \geq lk_{0}+1}\left(z,\xi,\nu\right), 
\\& \left(\nu+2\sqrt{-1}   \mathcal{L}\left(z,\xi\right)+Q_{k_{0}+1}\left(z,\xi,\nu\right)\right)^{l'}=\left(\nu+2\sqrt{-1}  \mathcal{L}\left(z,\xi\right)\right)^{l'}+\mbox{O}_{\mbox{wt} \geq l'k_{0}+1}\left(z,\xi,\nu\right). \end{split} 
\end{equation*}

It follows also that it does not exist any overlapping among the  following   terms
\begin{equation}\begin{split}&\quad\hspace{0.05 cm} f_{kl}z^{k}\mathcal{L}_{z}\left(z,\xi\right)\left(\nu+2\sqrt{-1}  \mathcal{L}\left(z,\xi\right)\right)^{l} ,\quad\mbox{for $k,l\in\mathbb{N}$ with $N=k+ k_{0}l-k_{0}+1$,}\\&   f_{k'l'}z^{k'}\mathcal{L}_{z}\left(z,\xi\right)\left(\nu+2\sqrt{-1}  \mathcal{L}\left(z,\xi\right)\right)^{l'} , \quad\mbox{for $k',l'\in\mathbb{N}$ with $N' = k' +  k_{0}l'-k_{0}+1$ and $2\leq N<N'$.} \end{split}\label{lolo2}
\end{equation}

In order to  compute (\ref{cete}), or equivalently (\ref{cete1}), 
we extract the sum of the  terms of weight $N\geq k_{0}+1$ from   (\ref{e3FFF}). We obtain
\begin{equation}\begin{split}&      \mathcal{L}_{z}\left(z,\xi\right)\displaystyle\sum_{k+k_{0}l=N-k_{0}+1}f_{kl}z^{k}\left(\nu+2\sqrt{-1}\mathcal{L}\left(z,\xi,\nu\right)\right)^{l} =\frac{\displaystyle\sum_{ k+k_{0}l=N}\left(g_{kl}z^{k}\left(\nu+2\sqrt{-1}\mathcal{L}\left(z,\xi\right)\right)^{l}\right)- \displaystyle\sum_{\tilde{k}+\tilde{l}=N}\tilde{g}_{\tilde{k}\tilde{l}}\xi^{\tilde{k}}\nu^{\tilde{l}}}{2\sqrt{-1}}     \\&   \quad\quad\quad\quad\quad\quad  +\\&    \quad\quad\quad\quad  \quad  
E_{N}\left(z,\xi,\nu\right) +\displaystyle\sum_{m+n+pk_{0}=N}\left(\varphi'_{mnp} - \varphi_{mnp}\right)z^{m}\xi^{n}\nu^{p} + \mathcal{L}_{\xi}\left(z,\xi\right)\displaystyle\sum_{\tilde{k}+\tilde{l}=N-k_{0}+1}\tilde{f}_{\tilde{k}\tilde{l}}\xi^{\tilde{k}}\nu^{\tilde{l}},\end{split}  \label{e3A}
\end{equation}
where we have  used the   notation  $E_{N}\left(z,\xi,\nu\right)= E_{N}\left(z,\xi,\nu;\left(f_{kl}\right)_{k+k_{0}l<N-k_{0}+1},\left(\tilde{f}_{\tilde{k}\tilde{l}}\right)_{\tilde{k}+ \tilde{l}<N-k_{0}+1};\left(g_{kl}\right)_{k+k_{0}l<N },\left(\tilde{g}_{\tilde{k}\tilde{l}}\right)_{\tilde{k}+ \tilde{l}<N}\right)$. 

It suffices by  (\ref{e3A})  to consider induction in respect to $N\geq k_{0}+1$, and  to focus on the coefficients of the following terms
\begin{equation}g_{kl}z^{k}\left(\nu+2\sqrt{-1}  \mathcal{L}\left(z,\xi\right)\right)^{l},\quad  \mathcal{L}_{z}\left(z,\xi\right)f_{k'l'}z^{k'}\left(\nu+2\sqrt{-1}  \mathcal{L}\left(z,\xi\right)\right)^{l'},\quad \xi^{\tilde{k}}\nu^{\tilde{l}}, \quad\mathcal{L}_{\xi}\left(z,\xi\right)\xi^{\tilde{k}'}\nu^{\tilde{l}'},\label{lolo3}
\end{equation}
for all $k,l,k',l',\tilde{k},\tilde{l},\tilde{k}',\tilde{l}'\in\mathbb{N}$ with $N=k+ k_{0}l$, $N = k' +  k_{0}l'-k_{0}+1$, $\tilde{k}'+\tilde{l}'=N-k_{0}+1$ and $\tilde{k}+\tilde{l}=N$. 
 
\subsection{\mbox{The case $s\neq 0$}} We study  terms in (\ref{e3FFF}) using  Model (\ref{mood}), when $s\neq 0$. In particular, we use the Model
 \begin{equation}w=\nu+2\sqrt{-1}  \nu^{s}  \mathcal{L}\left(z,\xi\right),\label{ref2}
\end{equation}

Procedures from Kol\`{a}\v{r}\cite{K1}, Zaitsev\cite{D1},\cite{D3} and Huang-Yin\cite{HY} are implemented in order  to consider   normalizations   in   (\ref{e3FFF}) when $s\neq 0$. More precisely, it is required to identify   interactions  containing undetermined terms of (\ref{cete11}) by linearising in (\ref{e3FFF}) like when  $s=0$. On the other hand, it is much more difficult to understand  interactions of terms, because there exist many  terms available in order to impose  normalizations like   the  terms derived from formal expansions involving   Model (\ref{ref2}).
  
Regardless of its non-triviality,  we want to make homogeneous   Model  (\ref{ref2}). It is introduced a system of pseudo-weights in this regard, because   it is  not possible to define a system of weights  like when $s=0$.   In particular, it is implemented the  strategy from \cite{V2}, but it is not clear how the system of pseudo-weights should be defined. In particular, we should have
\begin{equation}\mbox{wt}\left\{z\right\}=\mbox{wt}\left\{\xi\right\}=1,\quad\mbox{and then}\hspace{0.1 cm} \mbox{wt}\left\{\mathcal{L}\left(z,\xi\right) \right\}=k_{0}.\label{AA0}
\end{equation}
 
We should also have 
\begin{equation}\mbox{wt}\left\{\nu\right\}=k_{0},\quad\quad \mbox{wt}\left\{\nu^{s}\mathcal{L}\left(z,\xi\right) \right\}=k_{0}.\label{AA1} 
\end{equation}

In order to be satisfied as much of the axioms of the weight,  we should have
\begin{equation}\mbox{wt}\left\{\left(\nu+2\sqrt{-1}  \nu^{s}  \mathcal{L}\left(z,\xi\right)\right)^{n}\right\}=n\cdot \mbox{wt}\left\{\nu+2\sqrt{-1}  \nu^{s}  \mathcal{L}\left(z,\xi\right)\right\}.
\label{AA405} 
\end{equation}

Since  we want to make homogeneous   Model  (\ref{ref2}),    we should have
\begin{equation}\mbox{wt}\left\{\nu^{\alpha+s\beta}\left(  \mathcal{L}\left(z,\xi\right)\right)^{\beta}\right\}=nk_{0},\quad\mbox{for all $\alpha,\beta\in\mathbb{N}$ with $\alpha+\beta=n\in\mathbb{N}^{\star}.$  }\label{AA2} 
\end{equation}
 
It follows that  we should have
\begin{equation}\mbox{wt}\left\{z^{m}\xi^{p}  \nu^{\alpha+s\beta}\left(  \mathcal{L}\left(z,\xi\right)\right)^{\beta}\right\}=m+p+nk_{0},\quad\mbox{for all $\alpha,\beta,m,p\in\mathbb{N}$ with $\alpha+\beta=n\in\mathbb{N}^{\star}$,}\label{AA3} 
\end{equation}
because we should have
\begin{equation*}\mbox{wt}\left\{z^{m} \xi^{p}\nu^{\alpha+s\beta}\left(  \mathcal{L}\left(z,\xi\right)\right)^{\beta}\right\}=\mbox{wt}\left\{ \xi^{p}z^{m}\right\}+\mbox{wt}\left\{\nu^{\alpha+s\beta}\left(  \mathcal{L}\left(z,\xi\right)\right)^{\beta}\right\},\quad\mbox{for all $\alpha,\beta,m,p\in\mathbb{N}$ with $\alpha+\beta=n\in\mathbb{N}^{\star}$,} 
\end{equation*}
since it is desired by (\ref{AA405}) to be satisfied as much of the axioms 
of the weight. In particular, we should have
\begin{equation}\begin{split}& \quad\quad\quad\quad\quad\hspace{0.18 cm}\mbox{wt}\left\{   \xi^{\tilde{k}}\nu^{\tilde{l}+s}\right\}=k_{0}\left( s+\tilde{l}\right)+\tilde{k} ,\hspace{0.1 cm}  \quad\mbox{for all $\tilde{l},\tilde{k}\in \mathbb{N}$,}\\&\quad\quad\quad\quad\quad\quad \hspace{0.14 cm}\mbox{wt}\left\{z^{k}\nu^{s} \right\}=k+sk_{0},    \hspace{0.09 cm} \quad\quad\quad\quad\mbox{for all $k\in\mathbb{N}^{\star}$,}\\&\mbox{wt}\left\{\nu^{\alpha+s\beta+\gamma}\left(  \mathcal{L}\left(z,\xi\right)\right)^{\beta}\right\}=\left(n+\gamma\right)k_{0},\quad\quad\hspace{0.05 cm}\quad\mbox{for all $\alpha,\beta,\gamma\in\mathbb{N}$ with $\alpha+\beta=n\in\mathbb{N}^{\star}$,  } \end{split}
\label{AA4}
\end{equation}
because it is desired to have
\begin{equation*}\mbox{wt}\left\{\nu^{\gamma}\left(\nu+2\sqrt{-1}  \nu^{s}  \mathcal{L}\left(z,\xi\right)\right)^{n}\right\}=\mbox{wt}\left\{\left(\nu+2\sqrt{-1}  \nu^{s}  \mathcal{L}\left(z,\xi\right)\right)^{n}\right\}+\mbox{wt}\left\{\nu^{\gamma}\right\},\quad\mbox{for all $\alpha,\beta,\gamma\in\mathbb{N}$ with $\alpha+\beta=n\in\mathbb{N}^{\star}$.  }
\end{equation*}

Because   (\ref{AA2}) and (\ref{AA4}) should hold, we should have
\begin{equation} \begin{split}&  \mbox{wt}\left\{ \mathcal{L}_{\xi}\left(z,\xi\right) \xi^{\tilde{k}'}\nu^{\tilde{l}'+s}\right\}=\left(\tilde{l}'+1\right)k_{0}+\tilde{k}'-1,\quad\hspace{0.1 cm}\mbox{for all $\tilde{k}',\tilde{l}'\in \mathbb{N}$ with $\tilde{k}'\neq 0$.}  \\& \quad  \mbox{wt}\left\{ \mathcal{L}_{\xi}\left(z,\xi\right) \nu^{\tilde{l}'+s}\right\}=\left(\tilde{l}'+s\right)k_{0}+k_{0}-1, \quad\quad\mbox{for all $\tilde{l}'\in \mathbb{N}$.} \end{split}  \label{AA5} 
\end{equation}

It becomes clear by    (\ref{AA2}),(\ref{AA3}),(\ref{AA4}) and (\ref{AA5})   that we should have
\begin{equation}\begin{split}&\mbox{wt}\left\{\xi^{k} \mathcal{L}_{z}\left(z,\xi\right)\right\}=k+k_{0}-1,\quad\quad \mbox{wt}\left\{\xi^{k} \mathcal{L}_{\xi}\left(z,\xi\right)\right\}=k+k_{0}-1,\quad\quad\mbox{for all $k\in\mathbb{N}$,}\\& \mbox{wt}\left\{\nu^{k} \mathcal{L}_{z}\left(z,\xi\right)\right\}=kk_{0}+k_{0}-1,\quad  \mbox{wt}\left\{z^{k}\mathcal{L}_{z}\left(z,\xi\right)\right\}=k+k_{0}-1,\quad\quad\mbox{for all $k\in\mathbb{N}$,}\\& \mbox{wt}\left\{\nu^{k} \mathcal{L}_{\xi}\left(z,\xi\right)\right\}=kk_{0}+k_{0}-1,\quad  \mbox{wt}\left\{z^{k}\mathcal{L}_{\xi}\left(z,\xi\right)\right\}=k+k_{0}-1,\quad\quad\mbox{for all $k\in\mathbb{N}$.}\end{split}\label{AA10}  
\end{equation}
 
It becomes    motivating  to  estimate  the pseudo-weights   of other terms from (\ref{e3A}), but such computations depend on their writings  as products of lower order terms and their pseudo-weights should depend on  the  pseudo-weights of their lower order terms.  In particular,  in order to understand the most suitable value of  the following weights
\begin{equation}\begin{split}&\mbox{wt}\left\{ z^{k}\mathcal{L}_{z}\left(z,\xi\right)\nu^{s}\left(\nu+2\sqrt{-1} \nu^{s}  \mathcal{L}\left(z,\xi\right)\right)^{l}\right\},\hspace{0.1 cm}\quad\mbox{for all $k\in\mathbb{N}$,}\\&  \mbox{wt}\left\{ z^{k}\mathcal{L}_{z}\left(z,\xi\right)\left(\nu+2\sqrt{-1} \nu^{s}  \mathcal{L}\left(z,\xi\right)\right)^{l}\right\},\quad\quad\hspace{0.07 cm}\mbox{ for all $k,l\in\mathbb{N}$,}\\& \mbox{wt}\left\{ \nu^{s}\mathcal{L}_{z}\left(z,\xi\right)\left(\nu+2\sqrt{-1} \nu^{s}  \mathcal{L}\left(z,\xi\right)\right)^{l}\right\},\quad\quad\hspace{0.15 cm}\mbox{for all $l\in\mathbb{N}$,}\\& \mbox{wt}\left\{ z^{k}\nu^{s}\left(\nu+2\sqrt{-1} \nu^{s}  \mathcal{L}\left(z,\xi\right)\right)^{l}\right\},\quad\quad\quad\quad \quad\mbox{for all $k,l\in\mathbb{N}$,}\\&\mbox{wt}\left\{ z^{k}\nu^{s}\mathcal{L}_{z}\left(z,\xi\right)\right\},\quad\quad\quad\quad\quad\quad\quad\quad\quad \quad\quad \hspace{0.1 cm}  \mbox{for all $k\in\mathbb{N}$.}\label{pep11}\end{split}\end{equation}
 
Because  we can    define   the weight of a polynomial expression by different formulas,   the most suitable definition is not clear. We   move forward by \cite{V2} in order to define a  System of Pseudo-Weights. In particular, we define
\begin{equation}\begin{split}& \mbox{wt}\left\{\xi^{\alpha}\nu^{\beta}\right\}=\alpha+\beta k_{0},\quad\mbox{for all $\alpha,\beta\in\mathbb{N}$,} \\& \mbox{wt}\left\{z^{\alpha}\nu^{\beta}\right\}=\alpha+\beta k_{0},\quad\mbox{for all $\alpha,\beta\in\mathbb{N}$,}\\& \mbox{wt}\left\{z^{\alpha}\xi^{\beta}\right\}=\alpha+\beta,\quad\quad\hspace{0.035 cm}\mbox{for all $\alpha,\beta\in\mathbb{N}$.}\end{split}\label{pseu1}
\end{equation}
  
In order to make clear (\ref{AA2}) and (\ref{AA4}) and to extend (\ref{pseu1}), we define
\begin{equation} \mbox{wt}\left\{\nu^{N} z^{\alpha}\xi^{\beta}  \right\}=\left\{\begin{split}& N+\alpha+\beta,\quad\hspace{0.21 cm}\quad\mbox{for all $N, \alpha,\beta\in\mathbb{N}$ with
 $0<\alpha+\beta<k_{0}$,}\\& N-s+\alpha+\beta,\quad\mbox{for all $N, \alpha,\beta\in\mathbb{N}$ with $\alpha+\beta=k_{0}$ and  $\alpha\cdot \beta\neq 0$.}   \end{split}\right.\label{pseu3}
\end{equation}
 
We observe that the homogeneous polynomial $\mathcal{L}$ of degree  $k_{0}$, is defined  by  monomials  $ z^{a}\xi^{b}$  with  $a+b=k_{0}$ with $a,b\in\mathbb{N}^{\star}$.  We define 
\begin{equation} \mbox{wt}\left\{\nu^{N} z^{a\beta}\xi^{b\beta}  \right\}=\left(N-(s-1)\beta\right)k_{0},\quad\mbox{for all $N,\beta, a,b\in\mathbb{N}$ with $a+b=k_{0}$ and $\left(N-(s-1)\beta\right)k_{0}\geq0$.}\label{pseu4}
\end{equation}
 
It is required to define
\begin{equation}\begin{split}& \mbox{wt}\left\{\nu^{N} z^{a\beta+c}\xi^{b\beta}  \right\}=\left(N-(s-1)\beta\right)k_{0}+c,\quad\mbox{for all $N,c,\beta, a,b\in\mathbb{N}$ with $a+b=k_{0}$ and $\left(N-(s-1)\beta\right)k_{0}\geq0$,}\\& \mbox{wt}\left\{\nu^{N} z^{a\beta}\xi^{b\beta+c}  \right\}=\left(N-(s-1)\beta\right)k_{0}+c,\quad\mbox{for all $N,c,\beta, a,b\in\mathbb{N}$ with $a+b=k_{0}$ and $\left(N-(s-1)\beta\right)k_{0}\geq0$.}\end{split}\label{pseu5}
\end{equation}

Otherwise,  when $\left(N-(s-1)\beta\right)k_{0} \leq 0$, we define
\begin{equation} \mbox{wt}\left\{\nu^{N} z^{a\beta}\xi^{b\beta}  \right\}=\left(N-(s-1)\beta\right)k_{0}+\left(a+b\right)\left(\beta-\beta'\right), \label{pseu6}
\end{equation}
for all $N,  a,b\in\mathbb{N}$ with $a+b=k_{0}$ and $\beta'\in\mathbb{N}$ maximal such that $\left(N-(s-1)\beta'\right)k_{0}\geq0$.

The most suitable definition is attained when   the right-hand side of (\ref{pseu6}) is minimal, more precisely when $\beta'\in\mathbb{N}$ is maximal satisfying the above property, because there exist more evaluations  available. We denote them by
\begin{equation} \mbox{wt}_{N,a,b}\left\{\nu^{N} z^{a}\xi^{b}  \right\},\quad \mbox{where $N,a,b\in\mathbb{N}$ and $a,b\neq 0$.} \label{pseu7}
\end{equation}

It results that the best definition is   the following
\begin{equation}\mbox{wt}\left\{\nu^{N} z^{a}\xi^{b}  \right\}=\mbox{Min}\left(\mbox{wt}_{N,a,b}\left\{\nu^{N} z^{a}\xi^{b}  \right\}\right),  \quad \mbox{where $N,a,b\in\mathbb{N}$ and $a,b\neq 0$.} \label{pseu77}
\end{equation}  

It follows that Model   (\ref{ref2}) becomes pseudo-weighted-homogeneous in respect to the system of pseudo-weights from (\ref{pseu1}), (\ref{pseu3}),(\ref{pseu4}), (\ref{pseu5}),(\ref{pseu6}),(\ref{pseu7}),(\ref{pseu77}). It  remains to study the eventual overlappings of following terms
\begin{equation} \begin{split}&\quad
g_{kl}z^{k}\left(\nu+2\sqrt{-1} \nu^{s}  \mathcal{L}\left(z,\xi\right)+Q_{k_{0}+1}\left(z,\xi,\nu\right)\right)^{l},\quad\hspace{0.1 cm}\mbox{for $k,l \in\mathbb{N}$ with $N=k+k_{0}l$,}\\& g_{k'l'}z^{k'}\left(\nu+2\sqrt{-1} \nu^{s} \mathcal{L}\left(z,\xi\right)+Q_{k_{0}+1}\left(z,\xi,\nu\right)\right)^{l'},\quad\mbox{for $ k',l'\in\mathbb{N}$ with $N' = k' + k_{0}l'$ and $2 \leq  N < N' $.} \end{split} \label{loloC}
\end{equation}

These terms  do not overlap in (\ref{loloC}) because of their  different  pseudo-weights, in the light of the following pseudo-weighted expansions
\begin{equation}\begin{split}&  \hspace{0.1 cm}   \left(\nu+2\sqrt{-1} \nu^{s} \mathcal{L}\left(z,\xi\right)+Q_{k_{0}+1}\left(z,\xi,\nu\right)\right)^{l}= \left(\nu+2\sqrt{-1}\nu^{s}  \mathcal{L}\left(z,\xi\right)\right)^{l} + \mbox{O}_{\mbox{wt}\geq lk_{0}+1}\left(\nu,z,\xi\right),  \\&   \left(\nu+2\sqrt{-1} \nu^{s} \mathcal{L}\left(z,\xi\right)+Q_{k_{0}+1}\left(z,\xi,\nu\right)\right)^{l'}=\left(\nu+2\sqrt{-1}\nu^{s}\mathcal{L}\left(z,\xi\right)\right)^{l'} +\mbox{O}_{\mbox{wt}\geq lk_{0}+1}\left(\nu,z,\xi\right).  \end{split} \label{loloD}
\end{equation}

If follows that it does not exist any overlapping among the  following    terms
\begin{equation}\begin{split}& \quad \hspace{0.1 cm} f_{kl}z^{k}\mathcal{L}_{z}\left(z,\xi\right)\nu^{s}\left(\nu+2\sqrt{-1} \nu^{s} \mathcal{L}\left(z,\xi\right)\right)^{l} ,\quad \hspace{0.1 cm}\mbox{for $k,l \in\mathbb{N}$ with $N=k+ k_{0}l-k_{0}+1$,}\\&   f_{k'l'}z^{k'}\mathcal{L}_{z}\left(z,\xi\right)\nu^{s} \left(\nu+2\nu^{s}\sqrt{-1}  \mathcal{L}\left(z,\xi\right)\right)^{l'} ,\quad\hspace{0.12 cm} \mbox{for $k',l' \in\mathbb{N}$ with $N'=k'+ k_{0}l'-k_{0}+1$ and  $N <N'$.}\end{split}  \label{loloE}
\end{equation}
     
 In order to  compute (\ref{cete}), or equivalently (\ref{cete1}), 
we extract  terms of pseudo-weight $N\geq k_{0}+1$ from   (\ref{e3FFF}). We obtain    
 \begin{equation}\begin{split}&     \frac{\displaystyle\sum_{ k+k_{0}l=N}\left(g_{kl}z^{k}\left(\nu+2\sqrt{-1}\nu^{s}\mathcal{L}\left(z,\xi\right)\right)^{l}\right)- \displaystyle\sum_{\tilde{k}+\tilde{l}=N}\tilde{g}_{\tilde{k}\tilde{l}}\xi^{\tilde{k}}\nu^{\tilde{l}}
}{2\sqrt{-1}}= \mathcal{L}_{z}\left(z,\xi\right)\nu^{s} \cdot\displaystyle\sum_{k+k_{0}l=N-k_{0}+1\atop{k\neq 0}}f_{kl}z^{k}\left(\nu+2\sqrt{-1}\nu^{s}\mathcal{L}\left(z,\xi\right)\right)^{l}   +\mathcal{L}_{z}\left(z,\xi\right)\nu^{s} \cdot\\& \quad\quad\quad\quad\quad\hspace{0.15 cm}\quad\quad\quad\quad\quad\quad\quad\quad\quad\quad\quad\quad\quad\quad\quad\quad\quad\quad\quad\quad\quad\quad\quad\quad\quad f_{0\hspace{0.1 cm}\frac{N-s-k_{0}+1}{k_{0}}}\left(\nu+2\sqrt{-1}\nu^{s}\mathcal{L}\left(z,\xi\right)\right)^{\frac{N-s-k_{0}+1}{k_{0}}}   +  \mathcal{L}_{\xi}\left(z,\xi\right)\nu^{s}\cdot \\& \quad\quad\quad\quad\quad\quad\quad\quad\quad\quad \displaystyle\sum_{\tilde{k}+\tilde{l}=N-k_{0}+1\atop{\tilde{k}\neq 0}}\tilde{f}_{\tilde{k}\tilde{l}}\xi^{\tilde{k}}\nu^{\tilde{l}}+  \mathcal{L}_{\xi}\left(z,\xi\right)\nu^{s} \tilde{f}_{0\hspace{0.1 cm}N-k_{0}+1} \nu^{N-k_{0}+1}  +
\displaystyle\sum_{\mbox{of p.-w. $N$}}\left(\varphi'_{mnp} - \varphi_{mnp}\right)z^{m}\xi^{n}\nu^{p}+   E_{N} \left(z,\xi,\nu\right),\end{split}  \label{e3AB}
\end{equation}
where we have  used the   notation 
 $  E_{N}\left(z,\xi,\nu\right)= E_{N}\left(z,\xi,\nu;\left(f_{kl}\right)_{k+k_{0}l<N-k_{0}+1},\left(\tilde{f}_{\tilde{k}\tilde{l}}\right)_{\tilde{k}+ \tilde{l}<N-k_{0}+1};\left(g_{kl}\right)_{k+k_{0}l<N } \right)$.

It suffices  to consider by (\ref{e3AB}) induction,   with respect to $N\geq k_{0}+1$,
and to focus on the coefficients of the following terms
\begin{equation}\mathcal{L}_{\xi}\left(z,\xi\right) \xi^{\tilde{k}'}\nu^{\tilde{l}'+s},\quad \xi^{\tilde{k}}\nu^{\tilde{l}+s},\quad \mathcal{L}_{\xi}\left(z,\xi\right) \nu^{\tilde{l}'+s},\quad\nu^{N},\label{lolo3se}
\end{equation}
for all $\tilde{k},\tilde{l},\tilde{l}',\tilde{k}',\tilde{l}'\in\mathbb{N}$ with $\tilde{l}'=N-k_{0}-s+1$,   $\tilde{k}+\tilde{l}=N$ and $\tilde{k}'+\tilde{l}'=N-k_{0}+1$. 
 
The defined System of Pseudo-Weights for  $s\neq 0$ generalizes the defined System of Weights introduced for  $s=0$, because  (\ref{e3A}) reduces to (\ref{e3AB}) if  $s=0$, and  Model (\ref{ref2})  reduces to Model (\ref{ref1}) if  $s=0$. In particular, the terms  (\ref{lolo3se}) are   the terms (\ref{lolo3}) when $s=0$. They  can not be  used in order   to consider   classical Fischer Decompositions\cite{Sh},  but  
we make  computations in the both sides  from (\ref{e3A}) and (\ref{e3AB}) using:

  \section{Adapted Fischer-Decompositions}
 
 \subsection{Pseudo-Weighted Fischer Decompositions}  We consider  Fischer Decompositions like in  
(\ref{de}) and (\ref{degen}), but according to the previous system of pseudo-weights. The Fischer Differential Operator is appropriately defined.  In particular, we use    Pseudo-Weighted Fischer Decompositions    considering  $F\left(z,\xi,\nu\right)$ just a homogeneous polynomial and 
\begin{equation} \begin{split}& P_{1}\left(z,\xi,\nu\right)= z^{k}\left(\nu+2\sqrt{-1}  \nu^{s} \mathcal{L}\left(z,\xi\right)\right)^{l},\hspace{0.03 cm}\quad\quad\quad \quad\quad  \quad\quad\hspace{0.08 cm}\mbox{for all  $k,l \in\mathbb{N}$ such that  $k+k_{0}l=N \geq k_{0}+1,$} \\&     P_{2}\left(z,\xi,\nu\right)= \xi^{\tilde{k}}\nu^{\tilde{l}+s},\hspace{0.14 cm}\quad\quad\quad\quad\quad\quad\quad\quad\quad\quad\quad\quad\quad\quad\quad  \mbox{for all $ \tilde{k},\tilde{l} \in\mathbb{N}$   such that $\tilde{k}+\tilde{l}+s=N\geq k_{0}+1$,  }\\&    P_{3}\left(z,\xi,\nu\right)=z^{k'}\mathcal{L}_{z}\left(z,\xi\right)\nu^{s}\left(\nu+2\sqrt{-1}  \nu^{s} \mathcal{L}\left(z,\xi\right)\right)^{l'},\quad \quad\mbox{for all  $k',l'\in\mathbb{N}$  such that $k'+k_{0}l'+s=N-k_{0}+1$ and $k'\neq 0$, } \\&  P_{4}\left(z,\xi,\nu\right)=\mathcal{L}_{\xi}\left(z,\xi\right)  \xi^{\tilde{k}'}\nu^{\tilde{l}'+s},\quad\quad\quad\quad\quad\quad\quad\quad\quad\quad \hspace{0.03 cm}\quad\mbox{for all $\tilde{k}',\tilde{l}'\in\mathbb{N}$  such that  $\tilde{k}'+\tilde{l}'+s=N-k_{0}+1$ and $\tilde{k}'\neq 0$,}\\& P_{5}\left(z,\xi,\nu\right)=\mathcal{L}_{\xi}\left(z,\xi\right)\left(\nu+2\sqrt{-1}\nu^{s}\mathcal{L}\left(z,\xi\right)\right)^{\frac{N-s-k_{0}+1}{k_{0}}},\\& P_{6}\left(z,\xi,\nu\right)=\mathcal{L}_{\xi}\left(z,\xi\right) \nu^{N-k_{0}+s+1}. \end{split}\label{qer1VVV}
\end{equation}

\subsection{Weighted Fischer-Decompositions\cite{Sh}} We consider the Fischer Decompositions 
(\ref{de}) and (\ref{degen}) using the system of weights previously defined. In particular, the Differential Operator  (\ref{lele}) is  defined in respect to the system of weights previously defined:
$$\frac{\partial }{\partial \nu}\left(\nu\right):=\frac{\partial^{k_{0}} }{\partial \nu^{k_{0}} }\left(\nu\right)= k_{0} !,\quad\frac{\partial  }{\partial \nu}\left(\nu^{2}\right)=2^{k_{0}} k_{0} !\nu, $$
because it is differentiated towards to the weight of $\nu$, which is $k_{0}$.

Such Fischer Decompositions hold according to   classical Fischer Decompositions    recalled  in (\ref{de}) and (\ref{degen}). In particular, we use    Weighted Fischer Decompositions    considering $F\left(z,\xi,\nu\right)$ just a weighted-homogeneous polynomial and (\ref{qer1VVV}) holds and $s=0$.
 
The formal power series $\left(G_{1},G_{2},G_{3},G_{4},R\right)\left(z,\xi,\nu\right)$ are   homogeneous polynomials in $\left(z,\xi,\nu\right)$ in (\ref{degen}), and uniquely determined in (\ref{degen}) when (\ref{qer1VVV})  holds for $s=0$ or $s\neq 0$.  Then, we impose normalizations  in (\ref{e3FFF}) using
the weighted-homogeneous polynomials from (\ref{qer1VVV}) when  $s=0$ or when $s\neq 0$.     In particular, the Formal Equivalence (\ref{cete1}) is determined   as follows

\section{Computation of the Formal   Equivalence (\ref{cete1})} 
\subsection{\mbox{The Case $s\neq 0$}}   We impose   the following normalization condition
\begin{equation}\displaystyle\sum_{\mbox{pseudo-weight $ N$}}{\varphi'}_{mnp}{z}^{m}{\xi}^{n}\nu^{p}   \in S_{N},\quad\mbox{for all $N\geq k_{0}+1$,} \label{qer2V}
\end{equation}
by defining via (\ref{qer1VVV}) the intersection of Spaces  of Fischer Normalizations
\begin{equation*}S_{N}=S_{N}^{(1)}\cap S_{N}^{(2)}
\cap S_{N}^{(3)}\cap S_{N}^{(4)}\cap S_{N}^{(5)}\cap S_{N}^{(6)},\quad\mbox{for all $N\geq k_{0}+1$,}\end{equation*}
using the following-defined Spaces  of Fischer Normalizations
\begin{equation*}\begin{split}& S_{N}^{(1)}=    \displaystyle\bigcap_{   k+k_{0}l=N }\ker\left(z^{k}\left(\nu+2\sqrt{-1}\nu^{s}  \mathcal{L}\left(z,\xi\right)\right)^{l}\right)^{\star},
\quad\quad \quad\quad\quad\quad\quad\quad\quad\quad\quad\hspace{0.24 cm}  S_{N}^{(2)}= \displaystyle\bigcap_{    \tilde{k}+s+\tilde{l}=N  } \ker\left(\xi^{\tilde{k}}\nu^{\tilde{l}+s}\right)^{\star}, \\&  S_{N}^{(3)}= \displaystyle\bigcap_{k'+k_{0}l'+s=N-k_{0}+1 \atop{k'\neq 0} }    \ker\left(z^{k'}\nu^{s}\mathcal{L}_{z}\left(z,\xi\right)\left(\nu+2\sqrt{-1}\nu^{s}  \mathcal{L}\left(z,\xi\right)\right)^{l'}\right)^{\star}, \quad\quad S_{N}^{(4)}=  \displaystyle\bigcap_{ \tilde{k}'+s+\tilde{l}'=N-k_{0}+1\atop{\tilde{k}'\neq 0} }\ker\left(\mathcal{L}_{\xi}\left(z,\xi\right) \xi^{\tilde{k}'}\nu^{\tilde{l}'+s}\right)^{\star},     \end{split}
\end{equation*}
and respectively, the following-defined Spaces  of Fischer Normalizations
\begin{equation*}S_{N}^{(5)}=\ker\left(\mathcal{L}_{\xi}\left(z,\xi\right)\left(\nu+2\sqrt{-1}\nu^{s}\mathcal{L}\left(z,\xi\right)\right)^{\frac{N-s-k_{0}+1}{k_{0}}}\right)^{\star},\quad S_{N}^{(6)}=\ker\left(\mathcal{L}_{\xi}\left(z,\xi\right) \nu^{N-k_{0}+s+1}\right)^{\star}.
\end{equation*}

Then  (\ref{cete1}) is  determined by  induction   with respect to $N\geq k_{0}+1$.   We obtain:
\bp Let $\mathcal{M}$ and $\mathcal{M}'$  be defined as in  (\ref{gigi}) and (\ref{gigi1}) for $s\neq 0$. Then, it exists a  unique formal mapping defined as in (\ref{cete1})  transforming $\mathcal{M}$ into $\mathcal{M}'$, such that  (\ref{qer2V}) holds.  
\ep

\subsection{\mbox{The Case $s=0$}}   We impose   the following normalization condition
\begin{equation}\displaystyle\sum_{\mbox{weight $ N$}}{\varphi'}_{mnp}{z}^{m}{\xi}^{n}\nu^{p}   \in \tilde{S}_{N},\quad\mbox{for all $N\geq k_{0}+1$,} \label{pop}
\end{equation}
by defining via (\ref{qer1VVV}) the intersection of Spaces  of Fischer Normalizations
\begin{equation*}\tilde{S}_{N}=\tilde{S}_{N}^{(1)}\cap \tilde{S}_{N}^{(2)}
\cap \tilde{S}_{N}^{(3)}\cap \tilde{S}_{N}^{(4)} ,\quad\mbox{for all $N\geq k_{0}+1$,}\end{equation*}
using the following-defined Spaces  of Fischer Normalizations
\begin{equation*}\begin{split}& \tilde{S}_{N}^{(1)}=    \displaystyle\bigcap_{   k+k_{0}l=N }\ker\left(z^{k}\left(\nu+2\sqrt{-1}  \mathcal{L}\left(z,\xi\right)\right)^{l}\right)^{\star},
\quad\quad \quad\quad\quad\quad\quad\quad\quad\quad\quad  \tilde{S}_{N}^{(2)}= \displaystyle\bigcap_{    \tilde{k}+\tilde{l}=N  } \ker\left(\xi^{\tilde{k}}\nu^{\tilde{l}}\right)^{\star}, \\&  \tilde{S}_{N}^{(3)}= \displaystyle\bigcap_{k'+k_{0}l'=N-k_{0}+1  }    \ker\left(z^{k'}\nu^{s}\mathcal{L}_{z}\left(z,\xi\right)\left(\nu+2\sqrt{-1}   \mathcal{L}\left(z,\xi\right)\right)^{l'}\right)^{\star}, \quad\quad \tilde{S}_{N}^{(4)}=  \displaystyle\bigcap_{ \tilde{k}'+\tilde{l}'=N-k_{0}+1 }\ker\left(\mathcal{L}_{\xi}\left(z,\xi\right) \xi^{\tilde{k}'}\nu^{\tilde{l}'}\right)^{\star}.     \end{split}
\end{equation*}

Then  (\ref{cete1}) is  determined by  induction   with respect to $N\geq k_{0}+1$.   We obtain:
\bp Let $\mathcal{M}$ and $\mathcal{M}'$  be defined as in  (\ref{gigi}) and (\ref{gigi1}) for $s=0$. Then, it exists a  unique formal mapping defined as in (\ref{cete1})  transforming $\mathcal{M}$ into $\mathcal{M}'$, such that  (\ref{pop}) holds.  
\ep

The  Fischer Decomposition (\ref{qer2V}) is well-defined, regardless if $s=0$ or $s\neq 0$,
because the System of Pseudo-Weights, defined by (\ref{pseu1}),(\ref{pseu3}),(\ref{pseu4}),(\ref{pseu5}),(\ref{pseu6}),(\ref{pseu7}),(\ref{pseu77}), generalizes the System of Weights defined for $s=0$.
We move forward  to:
\subsection{Proof of Theorem \ref{t}} Denoting by $\mathcal{\tilde{M}}$ the  construction of normal form type, we use the following commutative diagram
\begin{equation}  
 \quad\begin{matrix}  \mathcal{H}_{1},\hspace{0.1 cm}\mathcal{H}_{2} \hspace{0.1 cm} :\mathcal{M} & \rightarrow&\hspace{0.1 cm}\mathcal{M}\\\quad\quad\quad\hspace{0.2 cm}\hspace{0.1 cm}\quad\Updownarrow& &\hspace{0.1 cm}\Updownarrow \\ \mathcal{H}_{1},\hspace{0.1 cm}\mathcal{H}_{2}\hspace{0.1 cm}:\tilde{\mathcal{M}}  &\rightarrow &\hspace{0.1 cm} \tilde{\mathcal{M}}  \end{matrix} ,\label{y56}
\end{equation}

 Any two Formal  Equivalences, which have  the same $1$-jet,   are actually identical because of the uniqueness of the Formal  Segre Equivalence sending $\mathcal{M}$ into $\mathcal{\tilde{M}}$. The  determination by $1$-jets follows from (\ref{y56}). 

\section{Analytic Equation}

We rewrite (\ref{cete1}) like   
$ \left(H(w,z),\tilde{H}\left(\nu,\xi\right)\right)=\left(F(w,z),G(w,z) ;F\left(\nu,\xi\right),G\left(\nu,\xi\right)\right)$ in order to rewrite (\ref{e3FFF}) with the notation
\begin{equation}\tilde{\varphi}'\left(w,z;\nu,\xi\right)= \displaystyle\sum_{m+n+p\geq k_{0}+1}\varphi'_{mnp} z^{m}  \xi^{n} \cdot  \left(\frac{ w+\nu}{2 }\right)^{p}.\label{brt}
\end{equation}

It follows that
\begin{equation}   \frac{ G(w,z)-\tilde{G}\left(\nu,\xi\right)}{2\sqrt{-1}} =\left( \frac{ G(w,z)+\tilde{G}\left(\nu,\xi\right)}{2 }\right)^{s}\cdot  \mathcal{L}\left( F(w,z),   \tilde{F}\left(\nu,\xi\right)\right)+\tilde{\varphi}'\left(G(w,z),F(w,z);\tilde{G}\left(\nu,\xi\right),\tilde{F}\left(\nu,\xi\right)\right), \label{e3FFFq}
\end{equation}
where $w=Q\left(z,\xi,\nu\right)$ such that (\ref{Qu}) holds, for 
\begin{equation}\tilde{\varphi}'\left(G(w,z),F(w,z);\tilde{G}\left(\nu,\xi\right),\tilde{F}\left(\nu,\xi\right)\right)= \displaystyle\sum_{m+n+p\geq k_{0}+1}\varphi'_{mnp} \left(   F(w,z)\right)^{m}\cdot\left(  \tilde{F}\left(\nu,\xi\right)  \right)^{n} \cdot  \left( \frac{ G(w,z)+\tilde{G}\left(\nu,\xi\right)}{2 }\right)^{p}.\label{brt1}
\end{equation}

It follows that
\begin{equation}\left.\Theta\left(w,z,\nu,\xi\right)\right|_{w=Q\left(z,\xi,\nu\right)}=0,\label{brt11}
\end{equation}
 such that (\ref{Qu}) holds,  where we have used the notation
\begin{equation}\Theta\left(w,z,\nu,\xi\right)=\frac{ G(w,z)-\tilde{G}\left(\nu,\xi\right)}{2\sqrt{-1}} -\left( \frac{ G(w,z)+\tilde{G}\left(\nu,\xi\right)}{2 }\right)^{s}\cdot  \mathcal{L}\left( F(w,z),   \tilde{F}\left(\nu,\xi\right)\right)-\tilde{\varphi}'\left(G(w,z),F(w,z);\tilde{G}\left(\nu,\xi\right),\tilde{F}\left(\nu,\xi\right)\right).\label{brt111}
\end{equation}

It is introduced the parametrization of maximum generic rank  
\begin{equation*}\rho\left(z,\xi,\nu\right)=\left(z,\xi;\nu+Q_{k_{0}}\left(z,\xi,\nu\right)+Q_{k_{0}+1}\left(z,\xi,\nu\right),\nu\right).\end{equation*}

 We apply  Proposition $6.2$ from Mir\cite{Mir1} in (\ref{brt11}),  because the right-hand side of (\ref{brt11}) is analytic. It follows  that $\Theta$ is analytic. In particular, we set $\nu=\xi=0$ 
 in (\ref{brt111}). It follows that
 \begin{equation}G(w,z)=\Theta\left(w,z;0,0\right). \label{1}
 \end{equation}
 
 Next, we set  $w=z=0$  in (\ref{brt111}). It follows that
  \begin{equation}\tilde{G}\left(\nu,\xi\right)=\Theta\left(0,0,\nu,\xi\right).  \label{2}
 \end{equation}

Then, (\ref{brt111}) is equivalent to
\begin{equation}\frac{ G(w,z)-\tilde{G}\left(\nu,\xi\right)}{2\sqrt{-1}}-\Theta\left(w,z,\nu,\xi\right)= \left( \frac{ G(w,z)+\tilde{G}\left(\nu,\xi\right)}{2 }\right)^{s}\cdot  \mathcal{L}\left( F(w,z),   \tilde{F}\left(\nu,\xi\right)\right)+\tilde{\varphi}'\left(G(w,z),F(w,z);\tilde{G}\left(\nu,\xi\right),\tilde{F}\left(\nu,\xi\right)\right), \label{vg}
\end{equation}
where its left-hand side is real-analytic, because (\ref{1}) and (\ref{2}) are real-analytic.

We make appropriate identifications of terms in the both sides of (\ref{vg}) using (\ref{brt1}). We obtain the convergence of $F(w,z)$ and $\tilde{F}\left(\nu,\xi\right)$. In particular, we obtain that (\ref{cete11}) is an analytic equivalence. In particular, (\ref{cete11}) is an analytic parametrization for (\ref{brt1}) and (\ref{vg}) is equivalent to

\begin{equation}\frac{ G(w,z)-\tilde{G}\left(\nu,\xi\right)}{2\sqrt{-1}}-\Theta\left(w,z,\nu,\xi\right)- \left( \frac{ G(w,z)+\tilde{G}\left(\nu,\xi\right)}{2 }\right)^{s}\cdot  \mathcal{L}\left( F(w,z),   \tilde{F}\left(\nu,\xi\right)\right)=\tilde{\varphi}'\left(G(w,z),F(w,z);\tilde{G}\left(\nu,\xi\right),\tilde{F}\left(\nu,\xi\right)\right). \label{vgg}
\end{equation}

 We apply  Proposition $6.2$ from Mir\cite{Mir1} in (\ref{vgg}) using (\ref{brt1}),  because its left-hand side   is real-analytic and its right-hand side is parametrized by an analytic equivalence. It follows that the formal power series (\ref{brt}) is convergent.

\subsection{Proof of Theorem \ref{t2}} It follows from the above computations related to (\ref{brt})-(\ref{vgg}).
\subsection{Proof of Corollary \ref{c2}}It follows like  above, but the formal power series (\ref{brt}) is assumed convergent.

\end{document}